\documentclass[11pt,a4paper]{amsart}
\pdfoutput=1
\usepackage[centering, margin=1in]{geometry}
\usepackage{amsmath,amsfonts,amsthm,amssymb,enumitem,graphicx,mathrsfs,mathtools}
\usepackage[stretch=10]{microtype}
\usepackage[unicode,breaklinks=true]{hyperref}
\usepackage{xcolor}
\definecolor{dark-red}{rgb}{0.4,0.15,0.15}
\definecolor{dark-blue}{rgb}{0.15,0.15,0.4}
\definecolor{medium-blue}{rgb}{0,0,0.5}
\hypersetup{
	pdftitle={Distributing Points on the Torus via Modular Inverses},
	pdfauthor={Peter Humphries},
	pdfnewwindow=true,
	colorlinks, linkcolor={dark-red},
	citecolor={dark-blue}, urlcolor={medium-blue}
}
\allowdisplaybreaks
\newcommand{\A}{\mathbb{A}}
\renewcommand{\C}{\mathbb{C}}
\newcommand{\CC}{\mathcal{C}}
\newcommand{\dee}{\partial}
\newcommand{\e}{\varepsilon}
\newcommand{\F}{\mathbb{F}}
\newcommand{\hf}{\mathbf{h}}
\newcommand{\R}{\mathbb{R}}
\newcommand{\RR}{\mathcal{R}}
\newcommand{\Ss}{\mathcal{S}}
\newcommand{\T}{\mathbb{T}}
\newcommand{\xf}{\mathbf{x}}
\newcommand{\Z}{\mathbb{Z}}
\DeclareMathOperator*{\vol}{vol}
\DeclareMathOperator*{\Var}{Var}
\numberwithin{equation}{section}
\newtheorem{theorem}[equation]{Theorem}
\theoremstyle{remark}
\newtheorem{remark}[equation]{Remark}
\begin{document}

\title[Distributing Points on the Torus via Modular Inverses]{Distributing Points on the Torus via Modular Inverses}

\author{Peter Humphries}

\address{Department of Mathematics, University of Virginia, Charlottesville, VA 22904, USA}

\email{\href{mailto:pclhumphries@gmail.com}{pclhumphries@gmail.com}}

\urladdr{\href{https://sites.google.com/view/peterhumphries/}{https://sites.google.com/view/peterhumphries/}}


\thanks{Research supported by the European Research Council grant agreement 670239.}

\begin{abstract}
We study various statistics regarding the distribution of the points
\[\left\{\left(\frac{d}{q},\frac{\overline{d}}{q}\right) \in \T^2 : d \in (\Z/q\Z)^{\times}\right\}\]
as $q$ tends to infinity. Due to nontrivial bounds for Kloosterman sums, it is known that these points equidistribute on the torus. We prove refinements of this result, including bounds for the discrepancy, small scale equidistribution, bounds for the covering exponent associated to these points, sparse equidistribution, and mixing.
\end{abstract}

\maketitle

\section{Introduction}

\subsection{Equidistribution}

For each positive integer $q$, consider the set of points in the torus $\T^2 \coloneqq (\R/\Z)^2$ given by
\begin{equation}
\label{eq:Sqdef}
S_q \coloneqq \left\{\left(\frac{d}{q},\frac{\overline{d}}{q}\right) \in \T^2 : d \in (\Z/q\Z)^{\times}\right\},
\end{equation}
where $\overline{d} \in (\Z/q\Z)^{\times}$ is the multiplicative inverse of $d$, so that $d\overline{d} \equiv 1 \pmod{q}$. There are $\varphi(q)$ such points; they are the image in $\T^2$ of the modulo $q$ hyperbola
\[\left\{\left(d,\overline{d}\right) \in ((\Z/q\Z)^{\times})^2 : d\overline{d} \equiv 1 \pmod{q}\right\}.\]
Associated to this set of points is the probability measure $\mu_q$ on $\T^2$ defined by
\begin{align*}
\mu_q(B) & \coloneqq \frac{1}{\varphi(q)} \#\left\{d \in (\Z/q\Z)^{\times} : \left(\frac{d}{q},\frac{\overline{d}}{q}\right)\in B\right\} \quad \text{for each Borel set $B \subset \T^2$,}	\\
\int_{\T^2} f(x) \, d\mu_q(x) & \coloneqq \frac{1}{\varphi(q)} \sum_{d \in (\Z/q\Z)^{\times}} f\left(\frac{d}{q},\frac{\overline{d}}{q}\right) \quad \text{for each measurable function $f : \T^2 \to \C$.}
\end{align*}
These measures equidistribute on $\T^2$ as $q \to \infty$ \cite{BK02,Zha96} (see also \cite{EL18}), so that
\begin{align}
\label{eq:equidistributionB}
\lim_{q \to \infty} \mu_q(B) & = \vol(B) \quad \text{for every continuity set $B \subset \T^2$,}	\\
\label{eq:equidistributionf}
\lim_{q \to \infty} \int_{\T^2} f(x) \, d\mu_q(x) & = \int_{\T^2} f(x) \, dx \quad \text{for every continuous function $f : \T^2 \to \C$.}
\end{align}
The proof of this fact is extremely short: it suffices show that for each fixed tuple $(m,n) \in \Z^2$,
\begin{equation}
\label{eq:equidistributionmn}
\lim_{q \to \infty} \int_{\T^2} e(mx_1 + nx_2) \, d\mu_q(x_1,x_2) = \int_{\T^2} e(mx_1 + nx_2) \, dx_1 \, dx_2 = \begin{dcases*}
1 & if $(m,n) = (0,0)$,	\\
0 & otherwise,
\end{dcases*}
\end{equation}
from which the Stone--Weierstrass theorem allows us to conclude \eqref{eq:equidistributionf}, at which point \eqref{eq:equidistributionB} follows from the Portmanteau theorem. Clearly \eqref{eq:equidistributionmn} holds for $(m,n) = (0,0)$. For $(m,n) \neq (0,0)$, the integral on the left-hand side of \eqref{eq:equidistributionmn} is simply $S(m,n;q)/\varphi(q)$, where
\[S(m,n;q) \coloneqq \sum_{d \in (\Z/q\Z)^{\times}} e\left(\frac{md + n\overline{d}}{q}\right)\]
denotes the Kloosterman sum, and so the Weil bound for Kloosterman sums (see \eqref{eq:Weil}) implies that this integral is $O_{m,n}(\tau(q) \sqrt{q}/\varphi(q))$ for $(m,n) \neq (0,0)$, from which equidistribution follows. (Note, however, that $S(m,n;q)/\varphi(q)$ can be much larger should $m$ and $n$ vary with $q$; in particular, if $m,n \equiv 0 \pmod{q}$, then this is \emph{equal} to $1$.)

In this paper, we study various refinements of this equidistribution result; we refer to the survey of Shparlinski \cite{Shp12} for further possible refinements and generalisations in different directions. Our emphasis is on quantifying in various ways how the measures $\mu_q$ behave like analogous measures associated to \emph{random} points. This is motivated by recent work of Bourgain, Rudnick, and Sarnak \cite{BRS17}, where analogous refinements of equidistribution are studied in the setting of lattice points on the sphere, namely statistics in the large $n$ limit of the projection onto the unit sphere $S^2 \subset \R^3$ of the set of
\[\left\{\left(x_1,x_2,x_3\right) \in \Z^3 : x_1^2 + x_2^2 + x_3^2 = n\right\}.\]

\subsection{Discrepancy}

Our first refinement is bounding the discrepancy of the measures $\mu_q$ as $q \to \infty$. The ball discrepancy is the quantity
\[D(\mu_q) \coloneqq \sup_{\substack{y \in \T^2 \\ 0 < R < \frac{1}{2}}} \left|\mu_q(B_R(y)) - \vol(B_R)\right|.\]
Here the supremum is over all injective geodesic balls $B_R(y)$ in $\T^2$. This is distinct from the box discrepancy $D^{\mathrm{box}}(\mu_q)$, where instead of taking a supremum over balls $B_R(y)$, one instead takes a supremum over all boxes $[a_1,b_1] \times [a_2,b_2]$ in $\T^2$; one can show that $D(\mu_q) \ll D^{\mathrm{box}}(\mu_q)^{1/2}$ via \cite[Theorem 2.1.6]{KN74}. It is natural to conjecture that $D(\mu_q) \ll_{\e} q^{-1/2 + \e}$, since this is the case for \emph{random} points. We make partial progress towards this conjecture, while also showing that for any fixed $\delta > 0$, the bound $|\mu_q(B_R(y)) - \vol(B_R)| \ll q^{-1/2 + \delta}$ is valid for \emph{almost all} centres $y \in \T^2$ of balls $B_R(y)$.

\begin{theorem}
\label{thm:discrepancy}
\hspace{1em}
\begin{enumerate}[leftmargin=*]
\item[\emph{(1)}] As $q \to \infty$, the ball discrepancy satisfies
\[D(\mu_q) \ll_{\e} q^{-\frac{1}{3} + \e}\]
for all $\e > 0$.
\item[\emph{(2)}] For $R < \frac{1}{2}$ and for any $\delta > 0$,
\[\vol\left(\left\{y \in \T^2 : \left|\mu_q(B_R(y)) - \vol(B_R)\right| > q^{-\frac{1}{2} + \delta}\right\}\right) \ll_{\e} \frac{q^{1 - 2\delta + \e}}{\varphi(q)}\]
for all $\e > 0$.
\end{enumerate}
\end{theorem}

This should be compared to the related work of Lubotzky, Phillips, and Sarnak on Hecke orbits of points on the sphere, where square-root cancellation of the spherical cap discrepancy is conjectured and the bound $O(q^{-1/3} (\log q)^{2/3})$ is proven \cite[Conjecture 2.4 and Theorem 2.5]{LPS86}. On the other hand, there exist number-theoretic situations where square-root cancellation is \emph{not} possible; Jung and Sardari have recently shown that the discrepancy associated to Hecke operators for modular forms of weight $k$ is $\Omega(k^{-1/3} (\log k)^{-2})$ \cite[Theorem 1.1]{JS20}.

\subsection{Small Scale Equidistribution}

We next consider small-scale equidistribution, namely the shrinking target problem in which one aims to show that $\mu_q(B_q)/\vol(B_q) \to 1$ for a sequence of sets $B_q$ whose volume shrinks as $q$ grows. By a pigeonhole-principle argument, we cannot always expect equidistribution in shrinking balls $B_R(y) \subset \T^2$ of radius $R$ for which $R = o(\varphi(q)^{-1/2})$. Moreover, there are specific regions where we \emph{never} find any points: if $d \in (\Z/q\Z)^{\times}$ with $2 \leq d \leq \sqrt{q}$, then $\overline{d} > \sqrt{q}$, so that $\mu_q(B_R(y)) = 0$ for $R = 1/2\sqrt{q}$ and $y = (1/2\sqrt{q},1/2\sqrt{q})$.

Nonetheless, it is natural to expect that equidistribution holds down to the optimal scale, so that for any fixed $\delta < 1/2$, we have that $\mu_q(B_R(y)) \sim \vol(B_R)$ for fixed $y \in \T^2$ and for all $q^{-\delta} < R < 1/2$. Such an asymptotic formula holds for \emph{random} points. The case $y = (0,0)$ would then imply a folklore conjecture on the existence of small modular inverses (see \cite{Gar06}), while the case $y = (1,0)$ would imply a slightly weaker form of a conjecture of Ford, Khan, Shparlinski, and Yankov on the maximal difference of modular inverses \cite[Conjecture 4.2]{FKSY05}.

Much like for the discrepancy, we are able to make partial progress towards these conjectures, as well as prove an optimal result for almost all centres of balls.

\begin{theorem}
\label{thm:smallscale}
\hspace{1em}
\begin{enumerate}[leftmargin=*]
\item[\emph{(1)}] Fix $y \in \T^2$ and $0 \leq \delta < \frac{1}{4}$. Then for $q^{-\delta} \leq R < \frac{1}{2}$, we have that
\[\mu_q(B_R(y)) = \vol(B_R) + O_{\e}\left(R^{\frac{2}{3}} q^{-\frac{1}{3} + \e}\right)\]
for all $\e > 0$.
\item[\emph{(2)}] Suppose that there exists some $\delta > 0$ such that $q^{-1/2 + \delta} \ll R < \frac{1}{2}$ as $q \to \infty$. Then for any fixed $\e > 0$,
\[\lim_{q \to \infty} \vol\left(\left\{y \in \T^2 : \left|\mu_q(B_R(y)) - \vol(B_R)\right| > \e \vol(B_R)\right\}\right) = 0.\]
\end{enumerate}
\end{theorem}

Related results for lattice points on the sphere have been proven by the author and Radziwi\l{}\l{} \cite[Theorem 1.5]{HR19}; see also \cite[Section 1.4]{BRS17}.

\subsection{Covering Exponents}

The \emph{covering radius} $\RR(P_n)$ of a set of points $P_n = \{x_1,\ldots,x_n\} \subset \T^2$ is the least $R > 0$ for which every point $y \in \T^2$ is within distance at most $R$ of some point $x_j$ in $P_n$. A packing argument implies that the covering radius of any set of $n$ points cannot be $o(1/\sqrt{n})$. The \emph{covering exponent} of a sequence $P$ of sets of points $P_n \subset \T^2$ is the quantity
\[K(P) \coloneqq -\liminf_{n \to \infty} \frac{\log n}{\log \vol(B_{\RR(P_n)})}.\]
Closely related to this is the \emph{average covering exponent} of a sequence $P = \{P_n\}$ of sets of points in $\T^2$. We let $\bar{\RR}(P_n,\delta)$ be the least $R > 0$ for which the measure of the set of points $y \in \T^2$ not within distance $R$ of a point in $P_n$ is at most $R^{-\delta}$. The average covering exponent of $P$ is
\[\bar{K}(P) \coloneqq - \lim_{\delta \to 0} \liminf_{n \to \infty} \frac{\log n}{\log \vol(B_{\bar{\RR}(P_n,\delta)})}.\]

\begin{theorem}
\label{thm:covering}
\hspace{1em}
\begin{enumerate}[leftmargin=*]
\item[\emph{(1)}] The covering exponent of the sequence of sets of points $S_q \subset \T^2$ in \eqref{eq:Sqdef} is at most $2$.
\item[\emph{(2)}] The average covering exponent of the sequence of sets of points $S_q \subset \T^2$ in \eqref{eq:Sqdef} is $1$.
\end{enumerate}
\end{theorem}

This should be compared to \cite[Section 1.4]{BRS17} and \cite[Section 1.1]{HR19} for analogous results for lattice points on the sphere $S^2$, and to \cite[Theorem 1.8 and Corollary 1.9]{Sar19a} and \cite[Corollary 1.6]{Sar19b} for lattice points on higher-dimensional spheres $S^d$ with $d \geq 3$.

\subsection{Variance Bounds and Asymptotics}

\hyperref[thm:discrepancy]{Theorems \ref*{thm:discrepancy} (2)} and \hyperref[thm:smallscale]{\ref*{thm:smallscale} (2)} are consequences of essentially sharp upper bounds for the variance
\begin{equation}
\label{eq:Var}
\Var(\mu_q;B_R) \coloneqq \int_{\T^2} \left(\mu_q(B_R(x)) - \vol(B_R)\right)^2 \, dx.
\end{equation}
It is natural to conjecture that for all $R \ll q^{-\delta}$ for some fixed $\delta > 0$, we have that
\[\Var(\mu_q;B_R) \sim \frac{\vol(B_R)}{\varphi(q)},\]
since such an asymptotic holds for \emph{random} points. The analogous statement for lattice points on the sphere is \cite[Conjecture 1.7]{BRS17}. A modification of this conjecture, replacing balls with annuli, was partially resolved by the author and Radziwi\l{}\l{} \cite[Theorem 1.3]{HR19}. We are able to resolve this conjecture provided that $q$ is prime and $R$ is sufficiently small, namely such that $R \sqrt{q} \log q \to 0$; note that generically $\mu_q(B_R(x)) = 0$ in this regime, so this should be thought of as the ``trivial'' regime. 

\begin{theorem}
\label{thm:variance}
\hspace{1em}
\begin{enumerate}[leftmargin=*]
\item[\emph{(1)}] 
For any $R < \frac{1}{2}$ and for all $\e > 0$, the variance \eqref{eq:Var} satisfies
\[\Var(\mu_q;B_R) \ll_{\e} \frac{\vol(B_R)}{\varphi(q)} q^{\e}.\]
\item[\emph{(2)}] Let $q$ be prime. For any $R < \frac{1}{2}$, the variance \eqref{eq:Var} satisfies
\[\Var(\mu_q;B_R) = \frac{\vol(B_R)}{\varphi(q)} + O\left(R^4 (\log q)^2\right).\]
\end{enumerate}
\end{theorem}

\subsection{Sparse Equidistribution}

We next consider the problem of sparse equidistribution, where we replace the measures $\mu_q$ with those associated to small subsets of $(\Z/q\Z)^{\times}$. Of course, equidistribution fails for subsets such as $\{a \in (\Z/q\Z)^{\times} : a \leq q/2\}$, since the corresponding measure is supported on $\{(x,y) \in \T^2 : x \leq 1/2\}$. For this reason, we restrict our study to subsets with \emph{algebraic} structure, namely cosets in $(\Z/q\Z)^{\times}$.

For each subgroup $H_q$ of $(\Z/q\Z)^{\times}$ and corresponding coset $aH_q \subset (\Z/q\Z)^{\times}$ (so that $a = 1$ corresponds to the subgroup itself), we define the probability measure $\mu_{aH_q}$ on $\T^2$ by
\begin{align*}
\mu_{aH_q}(B) & \coloneqq \frac{1}{\# H_q} \#\left\{d \in aH_q : \left(\frac{d}{q},\frac{\overline{d}}{q}\right)\in B\right\} \quad \text{for each Borel set $B \subset \T^2$,}	\\
\int_{\T^2} f(x) \, d\mu_{aH_q}(x) & \coloneqq \frac{1}{\# H_q} \sum_{d \in aH_q} f\left(\frac{d}{q},\frac{\overline{d}}{q}\right) \quad \text{for each measurable function $f : \T^2 \to \C$.}
\end{align*}
We prove the following.

\begin{theorem}
\label{thm:sparse}
Fix $\delta > 0$. For each positive cubefree integer $q$, pick a subgroup $H_q$ of $(\Z/q\Z)^{\times}$ and an associated coset $aH_q \subset (\Z/q\Z)^{\times}$ for which $\# H_q \gg q^{\frac{1}{2} + \delta}$. Then the probability measures $\mu_{aH_q}$ equidistribute on $\T^2$ as $q$ tends to infinity along cubefree integers. Furthermore, the same holds only under the assumption $\# H_q \gg q^{\delta}$ provided that $q$ tends to infinity along primes.
\end{theorem}

It is natural to conjecture that equidistribution holds under the weaker assumption $\# H_q \gg q^{\delta}$ for \emph{all} positive integers $q$, not just primes (though cf.~\hyperref[rem:failure]{Remark \ref*{rem:failure}}). This is the analogue of \cite[Conjecture 1]{MiVe06}, in which Michel and Venkatesh pose a similar conjecture for the equidistribution of subsets of Heegner points indexed by small subgroups of the class group of an imaginary quadratic field (see also \cite{HM06} and, more generally, \cite{Ven10}). Michel and Venkatesh note that in their setting, the generalised Lindel\"{o}f hypothesis implies a result analogous to \hyperref[thm:sparse]{Theorem \ref*{thm:sparse}}.

\subsection{Mixing}

Finally, we consider the problem of mixing. The group $(\Z/q\Z)^{\times}$ acts on the set $S_q$ in \eqref{eq:Sqdef} via
\[a \cdot \left(\frac{d}{q},\frac{\overline{d}}{q}\right) \coloneqq \left(\frac{ad}{q},\frac{\overline{ad}}{q}\right).\]
For each $a \in (\Z/q\Z)^{\times}$, we let
\begin{align*}
S_{q;a} & \coloneqq \left\{\left(\left(\frac{d}{q},\frac{\overline{d}}{q}\right), a \cdot \left(\frac{d}{q},\frac{\overline{d}}{q}\right)\right) \in \T^2 \times \T^2 : d \in (\Z/q\Z)^{\times}\right\}	\\
& = \left\{\left(\frac{d}{q},\frac{\overline{d}}{q},\frac{ad}{q},\frac{\overline{ad}}{q}\right) \in \T^4 : d \in (\Z/q\Z)^{\times}\right\}.
\end{align*}
We associate to this a probability measure $\mu_{q;a}$ on $\T^4 = \T^2 \times \T^2$ via
\begin{align*}
\mu_{q;a}(B) & \coloneqq \frac{1}{\varphi(q)} \#\left\{d \in (\Z/q\Z)^{\times} : \left(\frac{d}{q},\frac{\overline{d}}{q},\frac{ad}{q},\frac{\overline{ad}}{q}\right) \in B\right\} \quad \text{for each Borel set $B \subset \T^4$,}	\\
\int_{\T^4} f(x) \, d\mu_{q;a}(x) & \coloneqq \frac{1}{\varphi(q)} \sum_{d \in (\Z/q\Z)^{\times}} f\left(\frac{d}{q},\frac{\overline{d}}{q},\frac{ad}{q},\frac{\overline{ad}}{q}\right) \quad \text{for each measurable function $f : \T^4 \to \C$.}
\end{align*}
We are interested in the limiting behaviour of these probability measures.

\begin{theorem}
\label{thm:mixing}
The probability measures $\mu_{q;a}$ equidistribute on $\T^4$ as $q$ tends to infinity along primes if and only if $a$ and $q - a$ both tend to infinity with $q$.
\end{theorem}

This is the analogue of the mixing conjecture of Michel and Venkatesh on the joint equidistribution of Heegner points \cite{MiVe06}, which has been conditionally resolved by Khayutin \cite{Kha19}.

\section{Tools}

To begin, we let $k : \R^2 \times \R^2 \to \R$ be a point-pair invariant, so that $k(x + w,y + w) = k(x,y)$ for all $x,y,w \in \R^2$, which gives rise to a point-pair invariant $K : \T^2 \times \T^2 \to \R$ given by
\begin{equation}
\label{eq:Kxy}
K(x,y) \coloneqq \sum_{(m,n) \in \Z^2} k(x + (m,n),y).
\end{equation}
For $R > 0$, we take $k = k_R$ given by
\begin{equation}
\label{eq:kRxy}
k_R(x,y) \coloneqq \begin{dcases*}
1 & if $|x - y| \leq R$,	\\
0 & otherwise,
\end{dcases*}
\end{equation}
and let $K = K_R$ denote the associated point-pair invariant on $\T^2 \times \T^2$; if $R < 1/2$, then this is the indicator function of $B_R(y)$. We may calculate the Fourier coefficients of $K_R$ as follows:
\begin{align}
\notag
\widehat{K_R}((m,n),y) & \coloneqq \int_{\T^2} K_R(x,y) e(-(m,n) \cdot x) \, dx	\\
\notag
& = \int_{\R^2} k_R(x,y) e(-(m,n) \cdot x) \, dx	\\
\intertext{by \eqref{eq:Kxy},}
\notag
& = R^2 e(-(m,n) \cdot y) \int_{B_1(0)} e(-R(m,n) \cdot x) \, dx	\\
\intertext{upon making the change of variables $x \mapsto Rx + y$ and inserting \eqref{eq:kRxy},}
\notag
& = 2\pi R^2 e(-(m,n) \cdot y) \int_{0}^{1} J_0\left(2\pi R \sqrt{m^2 + n^2} r\right) r \, dr	\\
\intertext{by \cite[Chapter IV, Theorem 3.3]{SW71}, where $J_{\nu}(x)$ is a Bessel function,}
\label{eq:FourierKR}
& = \begin{dcases*}
\vol(B_R) & if $(m,n) = (0,0)$,	\\
\frac{R J_1\left(2\pi R \sqrt{m^2 + n^2}\right)}{\sqrt{m^2 + n^2}} e(-(m,n) \cdot y) & otherwise,
\end{dcases*}
\end{align}
by the fact that $J_0(0) = 1$ and \cite[6.561.5]{GR15}. The Fourier series for $K_R$ does not converge absolutely, since the Fourier coefficients are not of sufficiently rapid decay. To work around this issue, we consider
\[\widetilde{k}_{\rho}(x,y) \coloneqq \begin{dcases*}
\frac{1}{\vol(B_{\rho})} & if $|x - y| \leq \rho$,	\\
0 & otherwise.
\end{dcases*}\]
Then for $0 < \rho < R$, we define
\[k_{R,\rho}^{\pm}(x,y) \coloneqq k_{R \pm \rho} \ast \widetilde{k}_{\rho}(x,y) = \int_{\R^2} k_{R \pm \rho}(x,w) \widetilde{k}_{\rho}(w,y) \, dw.\]
It is readily checked that $k_{R,\rho}^{\pm}(x,y)$ are both nonnegative, continuous, pointwise linear in radial coordinates, bounded by $1$, and satisfy
\begin{align*}
k_{R,\rho}^{+}(x,y) & = \begin{dcases*}
1 & if $|x - y| \leq R$,	\\
0 & if $|x - y| > R + \rho$,
\end{dcases*}	\\
k_{R,\rho}^{-}(x,y) & = \begin{dcases*}
1 & if $|x - y| \leq R - \rho$,	\\
0 & if $|x - y| > R$.
\end{dcases*}
\end{align*}
Thus for all $x,y \in \T^2$, we have the pointwise inequalities
\begin{equation}
\label{eq:KRsqueeze}
K_{R,\rho}^{-}(x,y) \leq K_R(x,y) \leq K_{R,\rho}^{+}(x,y).
\end{equation}
Moreover, the Fourier coefficients are given by
\begin{equation}
\label{eq:FourierKRrhopm}
\widehat{K_{R,\rho}^{\pm}}((m,n),y) = \begin{dcases*}
\vol(B_{R \pm \rho}) & \hspace{-2.5cm} if $(m,n) = (0,0)$,	\\
\frac{(R \pm \rho)}{\pi \rho} \frac{J_1\left(2\pi (R \pm \rho) \sqrt{m^2 + n^2}\right) J_1\left(2\pi \rho \sqrt{m^2 + n^2}\right)}{m^2 + n^2} e(-(m,n) \cdot y) & \\
& \hspace{-2.5cm} otherwise;
\end{dcases*}
\end{equation}
in particular, the Fourier series for $K_{R,\rho}^{\pm}$ converges absolutely. This follows from the bound
\begin{equation}
\label{eq:J1bound}
J_1(x) \ll \min\left\{x,\frac{1}{\sqrt{x}}\right\}
\end{equation}
for $x \geq 0$ \cite[8.441.2 and 8.451.1]{GR15}.

\section{Proofs}

\begin{proof}[Proof of {\hyperref[thm:smallscale]{Theorem \ref*{thm:smallscale} (1)}}]
Since
\[\mu_q(B_R(y)) = \int_{\T^2} K_R(x,y) \, d\mu_q(x),\]
we have from \eqref{eq:KRsqueeze} that for any $0 < \rho < R$,
\begin{equation}
\label{eq:muqsqueeze}
\int_{\T^2} K_{R,\rho}^{-}(x,y) \, d\mu_q(x) \leq \mu_q(B_R(y)) \leq \int_{\T^2} K_{R,\rho}^{+}(x,y) \, d\mu_q(x).
\end{equation}
Via \eqref{eq:FourierKRrhopm}, we have the absolutely convergent spectral expansions
\begin{multline*}
\int_{\T^2} K_{R,\rho}^{\pm}(x,y) \, d\mu_q(x) = \pi (R \pm \rho)^2	\\
+ \frac{R \pm \rho}{\pi \rho \varphi(q)} \sum_{\substack{(m,n) \in \Z^2 \\ (m,n) \neq (0,0)}} \frac{J_1\left(2\pi (R \pm \rho\right) \sqrt{m^2 + n^2}) J_1\left(2\pi \rho \sqrt{m^2 + n^2}\right)}{m^2 + n^2} e(-(m,n) \cdot y) S(m,n;q).
\end{multline*}
We use the bounds \eqref{eq:J1bound} for the $J$-Bessel function and the Weil bound for Kloosterman sums,
\begin{equation}
\label{eq:Weil}
|S(m,n;q)| \leq \tau(q) \sqrt{(m,n,q) q},
\end{equation}
to bound the sum over $(m,n) \in \Z^2$ with $(m,n) \neq (0,0)$. Combining this with \eqref{eq:muqsqueeze}, we deduce that
\[\left|\mu_q(B_R(y)) - \vol(B_R)\right| \ll_{\e} R\rho + q^{-\frac{1}{2} + \e} + R^{\frac{1}{2}} \rho^{-\frac{1}{2}} q^{-\frac{1}{2} + \e}.\]
Upon taking $\rho = R^{-1/3} q^{-1/3}$ for $R > q^{-1/4}$ and $\rho = R/2$ for $R \leq q^{-1/4}$, we conclude that
\begin{equation}
\label{eq:discrepancybounds}
\left|\mu_q(B_R(y)) - \vol(B_R)\right| \ll_{\e} 
\begin{dcases*}
R^{\frac{2}{3}} q^{-\frac{1}{3} + \e} & for $R > q^{-\frac{1}{4}}$,	\\
q^{-\frac{1}{2} + \e} & for $R \leq q^{-\frac{1}{4}}$.
\end{dcases*}
\qedhere
\end{equation}
\end{proof}

\begin{remark}
When $q$ is prime, one can improve \eqref{eq:discrepancybounds} to
\[\left|\mu_q(B_R(y)) - \vol(B_R)\right| \ll
\begin{dcases*}
R^{\frac{2}{3}} q^{-\frac{1}{3}} & for $R > q^{-\frac{1}{4}}$,	\\
q^{-\frac{1}{2}} & for $R \leq q^{-\frac{1}{4}}$
\end{dcases*}\]
\end{remark}

\begin{proof}[Proof of {\hyperref[thm:discrepancy]{Theorem \ref*{thm:discrepancy} (1)}}]
This is an immediate consequence of \eqref{eq:discrepancybounds}.
\end{proof}

\begin{proof}[Proof of {\hyperref[thm:variance]{Theorem \ref*{thm:variance} (1)}}]
Via \eqref{eq:FourierKR}, the variance $\Var(\mu_q;B_R)$ has the absolutely convergent spectral expansion
\begin{equation}
\label{eq:Varspect}
\frac{R^2}{\varphi(q)^2} \sum_{\substack{(m,n) \in \Z^2 \\ (m,n) \neq (0,0)}} \frac{J_1\left(2\pi R\sqrt{m^2 + n^2}\right)^2}{m^2 + n^2} S(m,n;q)^2.
\end{equation}
We insert the bounds \eqref{eq:J1bound} for the $J$-Bessel function and the Weil bound for Kloosterman sums, \eqref{eq:Weil}, to see that this is bounded by a constant multiple of
\[\frac{R^4 \tau(q)^2 q}{\varphi(q)^2} \sum_{\substack{(m,n) \in \Z^2 \\ 0 < m^2 + n^2 \leq \frac{1}{R^2}}} (m,n,q) + \frac{R \tau(q)^2 q}{\varphi(q)^2} \sum_{\substack{(m,n) \in \Z^2 \\ m^2 + n^2 > \frac{1}{R^2}}} \frac{(m,n,q)}{(m^2 + n^2)^{3/2}},\]
which is easily seen to be $O_{\e}(R^2 q^{-1 + \e})$.
\end{proof}

\begin{remark}
Again, when $q$ is prime, one can improve \hyperref[thm:variance]{Theorem \ref*{thm:variance} (1)} to show that $\Var(\mu_q;B_R) \ll \vol(B_R) / \varphi(q)$.
\end{remark}

\begin{proof}[Proof of {\hyperref[thm:variance]{Theorem \ref*{thm:variance} (2)}}]
Upon opening up the Kloosterman sums in \eqref{eq:Varspect}, the variance is equal to the absolutely convergent spectral expansion
\[\frac{R^2}{\varphi(q)^2} \sum_{d_1,d_2 \in (\Z/q\Z)^{\times}} \sum_{\substack{(m,n) \in \Z^2 \\ (m,n) \neq (0,0)}} \frac{J_1\left(2\pi R\sqrt{m^2 + n^2}\right)^2}{m^2 + n^2} e\left(\frac{(d_1 - d_2)m}{q}\right) e\left(\frac{(\overline{d_1} - \overline{d_2})n}{q}\right).\]

The diagonal terms, namely those for which $d_1 = d_2$, contribute
\[\frac{R^2}{\varphi(q)} \sum_{\substack{(m,n) \in \Z^2 \\ (m,n) \neq (0,0)}} \frac{J_1\left(2\pi R\sqrt{m^2 + n^2}\right)^2}{m^2 + n^2} = \frac{\vol(B_R)}{\varphi(q)} - \frac{\vol(B_R)^2}{\varphi(q)}\]
by Parseval's identity, recalling \eqref{eq:FourierKR}.

The off-diagonal terms, namely those for which $d_1 \neq d_2$, break up into the sum of four separate terms dependent on the following conditions on $m$ and $n$:
\begin{enumerate}
\item $m \geq 1$ and $n \geq 0$,
\item $m \leq 0$ and $n \geq 1$,
\item $m \leq -1$ and $n \leq 0$,
\item $m \geq 0$ and $n \leq -1$.
\end{enumerate}
We observe that the third term is equal to the first term and the fourth term is equal to the second term by reindexing $(m,n)$ with $(-m,-n)$ and $(d_1,d_2)$ with $(-d_1,-d_2)$. Upon additionally reindexing $(m,n)$ with $(-n,m)$ and $(d_1,d_2)$ with $(\overline{d_1},\overline{d_2})$ for the second term, we deduce that the contribution from the off-diagonal terms is equal to the absolutely convergent expression
\begin{multline}
\label{eq:prepartsum}
\frac{4R^2}{\varphi(q)^2} \sum_{\substack{d_1,d_2 \in (\Z/q\Z)^{\times} \\ d_1 \neq d_2}} \lim_{X,Y \to \infty} \sum_{1 \leq m \leq X} e\left(\frac{(d_1 - d_2)m}{q}\right)	\\
\times \sum_{0 \leq n \leq Y} \cos\left(\frac{2\pi (\overline{d_1} - \overline{d_2})n}{q}\right) \frac{J_1\left(2\pi R\sqrt{m^2 + n^2}\right)^2}{m^2 + n^2},
\end{multline}
where we have used the fact that $e(x) + e(-x) = 2\cos(2\pi x)$.

We now use partial summation on both the sum over $1 \leq m \leq X$ and the sum over $0 \leq n \leq Y$ in \eqref{eq:prepartsum}, so that the inner double sum over $m$ and $n$ in \eqref{eq:prepartsum} is equal to
\begin{multline}
\label{eq:partsum}
\int_{0}^{Y} \int_{1}^{X} \sum_{1 \leq m \leq x} e\left(\frac{(d_1 - d_2)m}{q}\right) \sum_{0 \leq n \leq y} \cos\left(\frac{2\pi(\overline{d_1} - \overline{d_2})n}{q}\right) \frac{\dee^2}{\dee x \dee y} \frac{J_1\left(2\pi R\sqrt{x^2 + y^2}\right)^2}{x^2 + y^2} \, dx \, dy	\\
- \sum_{0 \leq n \leq Y} \cos\left(\frac{2\pi(\overline{d_1} - \overline{d_2})n}{q}\right) \int_{1}^{X} \sum_{1 \leq m \leq x} e\left(\frac{(d_1 - d_2)m}{q}\right) \frac{\dee}{\dee x} \frac{J_1\left(2\pi R\sqrt{x^2 + Y^2}\right)^2}{x^2 + Y^2} \, dx	\\
- \sum_{1 \leq m \leq X} e\left(\frac{(d_1 - d_2)m}{q}\right) \int_{0}^{Y} \sum_{0 \leq n \leq y} \cos\left(\frac{2\pi(\overline{d_1} - \overline{d_2})n}{q}\right) \frac{\dee}{\dee y} \frac{J_1\left(2\pi R\sqrt{X^2 + y^2}\right)^2}{X^2 + y^2} \, dy	\\
+ \sum_{1 \leq m \leq X} e\left(\frac{(d_1 - d_2)m}{q}\right) \sum_{0 \leq n \leq Y} \cos\left(\frac{2\pi(\overline{d_1} - \overline{d_2})n}{q}\right) \frac{J_1\left(2\pi R\sqrt{X^2 + Y^2}\right)^2}{X^2 + Y^2}.
\end{multline}
We shall show that the last three terms in \eqref{eq:partsum} converge to $0$ as $X,Y \to \infty$ and the first term is absolutely convergent as $X,Y \to \infty$; we shall then bound the limit of this first term.

By evaluating these geometric series, we have the bounds
\[\left|\sum_{1 \leq m \leq x} e\left(\frac{(d_1 - d_2)m}{q}\right)\right| \leq \frac{1}{\left|\sin \frac{\pi(d_1 - d_2)}{q}\right|}, \qquad \left|\sum_{0 \leq n \leq y} \cos\left(\frac{2\pi(\overline{d_1} - \overline{d_2})n}{q}\right)\right| \leq \frac{1}{\left|\sin \frac{\pi(\overline{d_1} - \overline{d_2})}{q}\right|}\]
independently of $x$ and $y$. Moreover, via \cite[8.440, 8.451.1, 8.471.1, and 8.471.2]{GR15},
\begin{align*}
\frac{\dee}{\dee x} \frac{J_1\left(2\pi R\sqrt{x^2 + y^2}\right)^2}{x^2 + y^2} & \ll_R \begin{dcases*}
x & for $x^2 + y^2 \ll_R 1$,	\\
\frac{x}{(x^2 + y^2)^2} & for $x^2 + y^2 \gg_R 1$,
\end{dcases*}	\\
\intertext{while}
\frac{J_1\left(2\pi R\sqrt{x^2 + y^2}\right)^2}{x^2 + y^2} & \ll_R \begin{dcases*}
1 & for $x^2 + y^2 \ll_R 1$,	\\
\frac{1}{(x^2 + y^2)^{3/2}} & for $x^2 + y^2 \gg_R 1$.
\end{dcases*}
\end{align*}
From this, the last three terms in \eqref{eq:partsum} tend to zero as $X,Y \to \infty$. For the first term in \eqref{eq:partsum}, we observe that
\begin{align*}
\int_{0}^{\infty} \int_{1}^{\infty} \left|\frac{\dee^2}{\dee x \dee y} \frac{J_1(2\pi R\sqrt{x^2 + y^2})^2}{x^2 + y^2}\right| \, dx \, dy & = (2\pi R)^2 \int_{0}^{\infty} \int_{2\pi R}^{\infty} \left|\frac{\dee^2}{\dee x \dee y} \frac{J_1(\sqrt{x^2 + y^2})^2}{x^2 + y^2}\right| \, dx \, dy	\\
& \ll R^2
\end{align*}
by first making the change of variables $x \mapsto x/2\pi R$ and $y \mapsto y/2\pi R$ and then once more using \cite[8.440, 8.451.1, 8.471.1, and 8.471.2]{GR15} in order to see that
\[\frac{\dee^2}{\dee x \dee y} \frac{J_1\left(\sqrt{x^2 + y^2}\right)^2}{x^2 + y^2} \ll \begin{dcases*}
xy & for $x^2 + y^2 \ll 1$,	\\
\frac{xy}{(x^2 + y^2)^{5/2}} & for $x^2 + y^2 \gg 1$.
\end{dcases*}\]

It follows that the contribution from the off-diagonal terms \eqref{eq:prepartsum} is bounded in absolute value by a constant multiple of
\[\frac{R^4}{\varphi(q)^2} \sum_{\substack{d_1,d_2 \in (\Z/q\Z)^{\times} \\ d_1 \neq d_2}} \frac{1}{\left|\sin \frac{\pi (d_1 - d_2)}{q} \sin \frac{\pi (\overline{d_1} - \overline{d_2})}{q}\right|} = \frac{R^4}{\varphi(q)^2} \sum_{c_1,c_2 \in (\Z/q\Z)^{\times}} \frac{\Ss(c_1,c_2)}{\left|\sin \frac{\pi c_1}{q} \sin \frac{\pi c_2}{q}\right|},\]
where
\[\Ss(c_1,c_2) \coloneqq \#\left\{(d_1,d_2) \in ((\Z/q\Z)^{\times})^2 : d_1 - d_2 \equiv c_1 \hspace{-.25cm} \pmod{q}, \ \overline{d_1} - \overline{d_2} \equiv c_2 \hspace{-.25cm} \pmod{q}\right\},\]
and we have used the fact that $q$ is prime to ensure that $c_1,c_2 \in (\Z/q\Z)^{\times}$. We have that
\[\Ss(c_1,c_2) = \#\left\{d_1 \in (\Z/q\Z)^{\times} : d_1^2 - c_1 d_1 + c_1\overline{c_2} \equiv 0 \hspace{-.25cm} \pmod{q}\right\},\]
and so $\Ss(c_1,c_2) \leq 2$ since the quadratic congruence $d_1^2 - c_1 d_1 + c_1\overline{c_2} \equiv 0 \pmod{q}$ has at most two solutions modulo a prime $q$ (see, for example, \cite[Lemma 9.6]{KL13}). We conclude that the contribution from the off-diagonal terms is bounded in absolute value by a constant multiple of
\[\frac{R^4}{\varphi(q)^2} \left(\sum_{c \in (\Z/q\Z)^{\times}} \frac{1}{\left|\sin \frac{\pi c}{q}\right|}\right)^2 \ll \frac{R^4}{\varphi(q)^2} q^2 (\log q)^2,\]
where the last inequality follows via the same method as the proof of the P\'{o}lya--Vinogradov inequality; see \cite[pp.~306--307]{MoVa07}.
\end{proof}

\begin{proof}[Proofs of {\hyperref[thm:discrepancy]{Theorems \ref*{thm:discrepancy} (2)}} and {\hyperref[thm:smallscale]{\ref*{thm:smallscale} (2)}}]
These follow from \hyperref[thm:variance]{Theorem \ref*{thm:variance}} via Chebyshev's inequality.
\end{proof}

\begin{proof}[Proofs of {\hyperref[thm:covering]{Theorems \ref*{thm:covering} (1) and (2)}}]
\hyperref[thm:covering]{Theorem \ref*{thm:covering} (1)} follows immediately from \hyperref[thm:smallscale]{Theorem \ref*{thm:smallscale} (1)}. Similarly, the upper bound for the average covering exponent in \hyperref[thm:covering]{Theorem \ref*{thm:covering} (2)} follows immediately from \hyperref[thm:smallscale]{Theorem \ref*{thm:smallscale} (2)}, while the lower bound is a simple consequence of a packing argument.
\end{proof}

\begin{proof}[Proof of {\hyperref[thm:sparse]{Theorem \ref*{thm:sparse}}}]
For each pair of integers $m,n \in \Z$, we have that
\begin{align}
\label{eq:aHqperiod}
\int_{\T^2} e(mx_1 + nx_2) \, d\mu_{aH_q}(x_1,x_2) & = \frac{1}{\# H_q} \sum_{d \in aH_q} e\left(\frac{md + n\overline{d}}{q}\right)	\\
\notag
& = \frac{1}{\varphi(q)} \sum_{\substack{\chi \hspace{-.25cm} \pmod{q} \\ \chi|_{H_q} = 1}} \overline{\chi}(a) \sum_{d \in (\Z/q\Z)^{\times}} \chi(d) e\left(\frac{md + n\overline{d}}{q}\right)
\end{align}
via character orthogonality. The sum over $d \in (\Z/q\Z)^{\times}$ is, by definition, the twisted Kloosterman sum $S_{\chi}(m,n;q)$. Since the number of characters $\chi$ modulo $q$ for which $\chi|_{H_q} = 1$ is $\varphi(q)/\# H_q$, the proof then follows from the Weil bound $S_{\chi}(m,n;q) \leq \tau(q) \sqrt{(m,n,q) q}$, which is known to hold for cubefree $q$ via \cite[Propositions 9.4, 9.7, 9.8, and Lemma 9.6]{KL13}.

When $q$ is prime, we instead note that since $H_q$ must be cyclic, we may write
\[\frac{1}{\# H_q} \sum_{d \in aH_q} e\left(\frac{md + n\overline{d}}{q}\right) = \frac{1}{q - 1} \sum_{x \in \F_q^{\times}} e\left(\frac{f(x)}{q}\right)\]
with $f(x) = a_1 x^{k_1} + a_2 x^{k_2}$ for $a_1 = am$, $a_2 = \overline{a}n$, $k_1 = (q - 1)/\# H_q$, and $k_2 = (q - 1)(\# H_q - 1)/\# H_q$. Since $q$ is a large prime and $m,n$ are fixed, we may assume without loss of generality that $m$ and $n$ are coprime to $q$. As
\[(k_1,q - 1) = (k_2,q - 1) = \frac{q - 1}{\# H_q}, \qquad (k_1 - k_2,q - 1) \leq 2\frac{q - 1}{\# H_q},\]
and $\# H_q \gg q^{\delta}$ by assumption, the conditions of \cite[Theorem 1]{Bou05} are met, which allows us to conclude that there exists $\delta' > 0$ such that
\[\frac{1}{\# H_q} \sum_{d \in aH_q} e\left(\frac{md + n\overline{d}}{q}\right) \ll q^{-\delta'}.\qedhere\]
\end{proof}

\begin{remark}
\label{rem:failure}
When $q$ is not cubefree, it is known that the Weil bound for $S_{\chi}(m,n;q)$ may fail; see \cite[Example 9.9]{KL13}. For this reason, it is conceivable that the extension of \hyperref[thm:sparse]{Theorem \ref*{thm:sparse}} to arbitrary $q$ is \emph{false} when $q$ is not cubefree.
\end{remark}

\begin{proof}[Proof of {\hyperref[thm:mixing]{Theorem \ref*{thm:mixing}}}]
It suffices to show that each for fixed $(m,n,m',n') \in \Z^4$,
\[S\left(m + am', n + \overline{a}n';q\right) = o_{m,n,m',n'}(q).\]
The Weil bound for Kloosterman sums \eqref{eq:Weil} implies that the left-hand side is bounded in absolute value by
\[2 \sqrt{(m + am', n + \overline{a}n',q) q}.\]
Since $q$ is a large prime and $m'$ and $n'$ are fixed, we may assume without loss of generality that $m'$ and $n'$ are coprime to $q$, in which case the Weil bound gives the desired result unless $a \equiv -m\overline{m'} \pmod{q}$ and $a \equiv -n\overline{n'} \pmod{q}$. However, since $m,m',n,n$ are fixed whereas $a$ and $q - a$ tend to infinity with $q$, these congruences cannot hold for sufficiently large $q$.

On the other hand, if $a$ does \emph{not} tend to infinity with $q$, then there exists a subsequence for which it is equal to a fixed positive integer $b$ by the Bolzano--Weierstrass theorem, and so
\[\int_{\T^4} e(-bx_1 - bx_2 + x_3 + x_4) \, d\mu_{q;a}(x_1,x_2,x_3,x_4) = \frac{S(-b + a, -b + a;q)}{\varphi(q)}\]
is equal to $1$ along this subsequence, which implies the failure of equidistribution, since
\[\int_{\T^4} e(-bx_1 - bx_2 + x_3 + x_4) \, dx_1 \, dx_2 \, dx_3 \, dx_4 = 0.\]
An analogous argument shows that equidistribution also fails if $q - a$ does not tend to infinity with $q$.
\end{proof}

\section{Generalisations}

Many of the results in this paper can be generalised to the setting studied by Granville, Shparlinski, and Zaharescu in \cite{GSZ05}. Following \cite[Section 2]{GSZ05}, we take an absolutely irreducible curve $\CC$ defined over $\F_p$ embedded in affine space $\A^r(\overline{\F}_p)$. We may naturally identify $\CC(\F_p)$ with a finite subset of $\T^r$ via the map $\F_p \to \T$ given by $x \mapsto x/p$. Let $\hf = (h_1,\ldots,h_s) : \CC \to \A^s(\overline{\F}_p)$ be a suitable rational map that is $L$-free along $\CC$ for some $L \geq 0$ in the sense of \cite[Section 2]{GSZ05}. Then \cite[Theorem 1]{GSZ05} states that the probability measures on $\T^s$ given by
\begin{align*}
\mu_{\CC,\hf}(B) & \coloneqq \frac{1}{\# \CC(\F_p)} \#\left\{\xf \in \CC(\F_p) : \hf(\xf) \in B\right\} \quad \text{for each Borel set $B \subset \T^s$,}	\\
\int_{\T^s} f(x) \, d\mu_{\CC,\hf}(x) & \coloneqq \frac{1}{\# \CC(\F_p)} \sum_{\xf \in \CC(\F_p)} f(\hf(\xf)) \quad \text{for each measurable function $f : \T^s \to \C$}
\end{align*}
equidistribute on $\T^s$ as $p$ tends to infinity provided that $L$ tends to infinity with $p$.

The case of $\CC$ being the curve $y = x$ and $\hf(x) = (x,x^{-1})$ corresponds to the equidistribution of the set $S_q \subset \T^2$ in \eqref{eq:Sqdef} when $q = p$ is prime. The key tool behind \cite[Theorem 1]{GSZ05} is a bound for exponential sums due to Bombieri \cite[(8)]{GSZ05}, which includes the Weil bound for Kloosterman sums of prime level as a special case. Since the key tools for several of the results in this paper are the Weil bound together with bounds for the Fourier coefficients of the indicator function of a ball in $\T^2$, which of course can be generalised to $\T^s$ with $s \geq 3$, it follows that many of the results in this paper can be generalised to the setting studied by Granville, Shparlinski, and Zaharescu in \cite{GSZ05}. Notably, the discrepancy bound proven in \hyperref[thm:discrepancy]{Theorem \ref*{thm:discrepancy} (1)} is proven in further generality in \cite[Lemma 3]{GSZ05}, albeit for the box discrepancy instead of the ball discrepancy.

It is less clear if these methods generalise to curves over $\Z/q\Z$ with $q$ nonsquarefree, for then the method of Bombieri \cite[(8)]{GSZ05} used to bound exponential sums is no longer valid, and instead one must use more elementary methods (see, for example, \cite[Lemmata 12.2 and 12.3]{IK04}). Notably, it is not necessarily the case that one can expect bounds in the nonsquarefree setting that are as strong as in the squarefree setting; cf.~\hyperref[rem:failure]{Remark \ref*{rem:failure}}.

\phantomsection
\addcontentsline{toc}{section}{Acknowledgements}
\hypersetup{bookmarksdepth=-1}

\subsection*{Acknowledgements}

Thanks are owed to Simon Rydin Myerson, to Andrew Granville for helpful discussions regarding \cite{GSZ05}, to Igor Shparlinski for alerting the author to \cite{Bou05} and its application in \hyperref[thm:sparse]{Theorem \ref*{thm:sparse}}, and especially to the anonymous referee for pointing out a serious error in an earlier version of this paper.

\hypersetup{bookmarksdepth}

\end{document}